\documentclass[a4paper,11pt]{article}
\usepackage[english]{babel}
\usepackage[psamsfonts]{amssymb}
\usepackage{cmmib57}
\usepackage{exscale}
\usepackage{amsmath}
\usepackage{amscd}
\usepackage{latexsym}
\usepackage{graphicx}
\usepackage{amsmath}
\newtheorem{theorem}{Theorem}

\newtheorem{corollary}[theorem]{Corollary}

\newtheorem{definition}[theorem]{Definition}

\newtheorem{proposition}[theorem]{Proposition}
\newtheorem{remark}[theorem]{Remark}

\newcommand{\PP}{\mathbb{P}}
\newcommand{\EE}{\mathbb{E}}

\newcommand{\no}{\noindent}
\newcommand{\ds}{\displaystyle}

\def\squarebox#1{\hbox to #1{\hfill\vbox to #1{\vfill}}}
\newcommand{\qed}{\hspace*{\fill}
\vbox{\hrule\hbox{\vrule\squarebox{.45em}\vrule}\hrule}\smallskip}

\begin{document}
\title{A note on the coupon - collector's problem
with multiple arrivals and the random sampling}
\author{
Marco Ferrante\footnote{corresponding author} \ and Nadia Frigo
\\
Dipartimento di Matematica \\
Universit\`a degli Studi di Padova \\
via Trieste, 63\\
35121 Padova, Italy \\
e-mail: ferrante@math.unipd.it  \ and \ nadia.frigo@gmail.com}
\maketitle
\begin{abstract}
\no
In this note we evaluate the expected
waiting time to complete a collection of coupons, in the
case of coupons which arrives in groups of constant size,
independently and with unequal probabilities.
As an application we will be able to determine the expected
number of samples of dimension g that we have to draw
independently in order to observe all the types of individuals
in a given population.
\end{abstract}

\vspace{3truecm}

%\vspace{1truecm}
%  {\bf Short title:} Coupon collector

\section{Introduction}
The coupon-collector's problem is a classical problem in combinatorial probability.
The description of the basic problem is easy:
consider one person that collects coupons and
assume that there is a finite number, say $m$, of different types of coupons.
These items arrive one by one in sequence,
with the type of the successive items being independent random variables
that are each equal to $k$ with probability $p_k$.
It is immediate to see how this description can be adapted
to the general problem to draw independent samples
from a given, finite distribution.

In the coupons-collector's problem,
one is usually interested in answering the following questions:
which is the probability to complete the collection
(or a given subset of the collection)
after the arrival of exactly $n$ coupons ($n \ge m$)?
which is the expected number of coupons that we need to complete the collection (or to complete
a given subset of the collection)? how these probabilities and expected values change if we
assume that the coupons arrive in groups of constant size?

The first results, due to De Moivre, Laplace and
Euler (see \cite{MR1082197} for a comprehensive introduction on this topic),
deal with the case of
constant probabilities $p_k\equiv \frac{1}{m}$,
while the first results on the unequal case have to be ascribed to Von Schelling
(see \cite{MR0061772}).
Many other studies have been carry out on this classical problem ever since
(see e.g. \cite{MR0268929}, \cite{MR959649} and \cite{MR1978107}).

The aim of this note is to evaluate
the expected number of coupons that one needs to collect
in order to complete the collection,
in the case of unequal probabilities and multiple arrival (i.e. the case
in which the coupons arrives in groups of constant size).
To the best of our knowledge
this result is new and
in the case of uniform probabilities
we derive the expression present in the literature
(see e.g. Stadje \cite{MR1082197})
in a much easier way.
Furthermore, we will apply this computation to the problem to
sample without replacement $g$ individuals from a population
composed by $m$ types of individuals, present in different proportions,
obtaining an explicit computation of the expected
number of independent samples of size g
that we have to draw in order to observe all the types of individuals in the
population, which could be of interest for other applied problems.

\section{The single arrival case}
In order to solve the problem for the multiple arrival setting,
we shall start by the easier single arrival case
and we shall see in the next section how to extend
this result to that case.

Let us fix the notation. We shall denote by $\{1, \ldots, m\}$ the
different types of items which form the collection
% (where $m \in \Z $ will also denote
% the total number of different types present in the collection)
Let us assume that the items are purchased one by one in sequence,
with the type of the successive items being independent random variables
that are each equal to $k$ with probability $p_k$.
Since we are interested here in the number of items one needs to collect
to complete the collection, let us define the following set of random variables:
$X_1$ will denote the (random) number of items that we need to collect to have the
first coupon of our collection (which is trivially equal to 1),
$X_2$ will be the number of
additional items that we need to collect to obtain the second
different coupon in our collection and
so on let us define, for every $i\le m$, by $X_i$ the number of
items that we need to collect to pass form the $i-1$-th to the $i$-th different coupon in the collection.
From this classical description (see e.g. Rosen \cite{MR0268929}),
we obtain that the random number
of coupons that we need to complete the collection
is equal to $X=X_1+\ldots + X_m$ and that $\PP[X<+\infty]=1$.

In the case of constant probabilities, i.e.
$p_k \equiv 1/m$ for any $k\in\{1,\ldots ,m\}$,
it is immediate to see that
the random variable $X_i$, for $i\in \{2, \ldots, m\}$, has a geometric law with
parameter $(m-i)/m$. The expected number of coupons that we
need in order to complete the collection is therefore given by the well-known formula
\begin{equation}
\label{1}
\EE[X] = m \sum_{i=1}^{m} \frac{1}{i}
\quad.
\end{equation}

When the probabilities $p_k$ are unequal,
one can look at the problem from a slightly different angle.
Let us define the following set of random variables:
$Y_1$ will denote the (random) number of items that we need to collect to obtain the
first coupon of type $1$, $Y_2$ the number of items that we need to collect
to get the first coupon of type $2$, and so on for the others coupons.
In this setting, the waiting time to complete the collection is given by
the random variable $Y=\max(Y_1, \ldots, Y_m)$.
In order to compute its expected value, one
can use the Maximum-Minimums identity (see \cite{MR732623}, p.345), obtaining
\begin{equation}
\label{1.5}
\begin{array}{rl}
\EE[Y] =
&
{\ds \sum_i \EE[Y_i] - \sum_{i<j} \EE[\min(Y_i,Y_j)] +
\sum_{i<j<k} \EE[\min(Y_i,Y_j,Y_k)] + \ldots}
\\
&
{\ds \ldots + (-1)^{m+1} \EE[\min(Y_1,Y_2,\ldots, Y_m)]
\quad.}
\end{array}
\end{equation}
Since the random variables $\min(Y_{i_1},Y_{i_2},\ldots, Y_{i_k})$
have a geometric law with parameter $ p_{i_1}+p_{i_2}+\ldots +p_{i_k}$, we
get the formula
\begin{equation}
\label{2}
\EE[Y] = \sum_i \frac{1}{p_i} - \sum_{i<j}
\frac{1}{p_i+p_j} +
\sum_{i<j<k}
\frac{1}{p_i+p_j+p_k} +
\ldots + (-1)^{m+1} \frac{1}{p_1+\ldots + p_m}
.
\end{equation}
The problem described above can be rephrased as follows:
let us consider a finite distribution and let us evaluate the
expected number of independent sample that we
have to draw in order to observe all the records.
The quantity (\ref{2}) is clearly this value.
It is interesting to note that if we would like to evaluate the
expected number of independent samples that we
have to draw in order to observe a fixed number $k$ of records,
with $k\le n$,
the present approach is no more suitable, but
we have to reconsider
the problem from a slightly different point of view
(see \cite{FerranteFrigo} for the details).

\section{The multiple arrival case}
Let us now consider the case of
coupons which arrives in groups of constant size $g$, where $1<g<m$,
with the types of the items in any group of
coupons being independent random variables.
A natural requirement in this contest is that
each group does not contain more than
one coupon of any type.
With this assumption, the total number of groups will be
${m \choose g}$ and each group $A$ can be identified with a vector
$(a_1,\ldots, a_g) \in \{1, \ldots, m\}^g$
with $a_i < a_{i+1}$ for $i=1, \ldots, g-1$.
Removing this assumption, we have to consider all the possible
$m^g$ groups of coupons that we can obtain. In this case
we will describe the groups of coupons as
$(a_1,\ldots, a_g) \in \{1, \ldots, m\}^g$ and we will
see at the end of this chapter how this problem
applies to the case of sampling form
a given population.

Let us first consider
the case in which each group does not contain more than
one coupon of any type.
We can order the groups according to the
lexicographical order (i.e.
$A=(a_1,\ldots, a_g) < B=(b_1,\ldots, b_g)$ if there exists
$i\in\{1, \ldots, g-1\}$ such that
$a_s = b_s$ for $s<i$ and $a_i < b_i$).
\begin{definition}
\label{def1}
We shall denote
by $q_i, i\in \{1, \ldots, {m \choose g}\}$ the probability
to purchase (at any given time)
the $i$-th group of coupons, accordingly to the lexicographical order.
Moreover, given $k\in \{1,\ldots,m-g\}$, we shall denote by
$q(i_1,\ldots,i_k)$ the probability to purchase a group of coupons which
does not contain any of the coupons $i_1,\ldots,i_k$.
\end{definition}
\begin{remark}
\label{remark2}
In order to compute the probabilities $q(i_1,\ldots,i_k)$'s, one
can proceed as follows: by the defined ordering, it holds that
\[
q(1)=
\sum_{i={m-1 \choose g-1}+1}^{{m \choose g}} q_i \ \ , \ \
q(1,2)=
\sum_{i={m-1 \choose g-1}+{m-2 \choose g-1}+1}^{{m \choose g}} q_i \ , \
\]
and in general
\[
q(1,2,\ldots,k)=
\left\{
\begin{array}{ll}
{\ds
\sum_{i={m-1 \choose g-1}+\cdots+{m-k \choose g-1}+1}
^{{m \choose g}} q_i \ , \ }
&
\mbox{if} \ k\le m-g
\\
0
&
\mbox{otherwise} \ .
\end{array}
\right.
\]
For any permutation $(i_1,\ldots,i_m)$ of $(1,\ldots,m)$, one first
reorders the $q_i$'s according to the lexicographical order of this new alphabet
and then compute
\[
q(i_1,i_2,\ldots,i_k)=
\left\{
\begin{array}{ll}
{\ds
\sum_{i={m-1 \choose g-1}+\cdots+{m-k \choose g-1}+1}
^{{m \choose g}} q_i \ , \ }
&
\mbox{if} \ k\le m-g
\\
0
&
\mbox{otherwise} \ .
\end{array}
\right.
\]
\end{remark}
\begin{remark}
\label{remark3}
There are many conceivable choices
for the unequal probabilities $q_i$'s.
For example, we can assume that one forms the groups following
the strategy of the draft lottery in the American professional sports, where
different proportion of the different coupons are put together and
we choose at random in sequence the coupons, discarding the eventually duplicates,
up to obtaining a group of $k$ coupons.
Or, more simply, we
can assume that
the $i$-th coupon will arrive with probability $p_i$ and that the
probability of any group is proportional to the product of the
probabilities of the single coupons contained.
\end{remark}

In order to evaluate the expected number of groups needed in order
to complete the collection, we shall use
the approach of the single arrival case.
Let us start by considering the case of uniform probabilities,
i.e. $q_i=\frac{1}{{m \choose g}}$ for any $i$.
Let us define the following set of random variables:
\[
V_i = \{\mbox{number of groups to purchase
to obtain the first coupon of type $i$}\}
\]
These random variables have a geometric law with
parameter
\[
1-\frac{{m-1 \choose g}}{{m \choose g}} \ .
\]
The random variables $\min(V_i,V_j)$ have their selves a geometric law
with parameter
\[
1-\frac{{m-2 \choose g}}{{m \choose g}}
\]
and so on up to the random variables $\min(V_{i_1},\ldots,V_{i_{m-g}})$,
which have a geometric law
with parameter
\[
1-\frac{1}{{m \choose g}} \ .
\]
The minimum of more random variables, i.e. $\min(V_{i_1},\ldots,V_{i_{k}})$ for $k>m-g+1$,
will be equal to the constant random variable $1$.

Applying the Maximum-Minimums principle, we shall obtain
that the expected number of groups of coupons that we need to complete
the collection is equal to
\[
\EE[\max(V_1, \ldots, V_m)]
=
\sum_{1\le i\le m} \EE[V_i]
- \sum_{1\le i<j\le m} \EE[\min(V_i,V_j)]
+ \ldots
\]
\[
+ (-1)^{m-g+1}
\sum_{1\le i_1<i_2<\cdots<i_{m-g}\le m}
\EE[V_{i_1},\ldots,V_{i_{m-g+1}}]
+
\]
\[
+(-1)^{m-g+2}
\sum_{1\le i_1<i_2<\cdots <i_{m-g+1}\le m} \ 1
+ \ldots + (-1)^{m+1}
\]
\[
=
{m \choose 1}
\frac{1}
{1-\frac{{m-1 \choose g}}{{m \choose g}}}
-
{m \choose 2}
\frac{1}
{1-\frac{{m-2 \choose g}}{{m \choose g}}}
+
{m \choose 3}
\frac{1}
{1-\frac{{m-3 \choose g}}{{m \choose g}}}
+ \ldots
\]
\[
\ldots +
(-1)^{m-g+1}
{m \choose m-g}
\frac{1}
{1-\frac{1}{{m \choose g}}}
+
\sum_{1\le k\le g}(-1)^{m-g+k+1}
{m \choose m-g+k} \ .
\]
This result, even if not obtained with this computation,
is known (see e.g. Stadje \cite{MR1082197}, p.872).

In the unequal case, we are able to generalize the previous result
as follows:
\begin{proposition}
\label{3.1}
The expected number of groups of coupons that we need to complete
the collection, in the case of unequal probabilities $q_i$,
is equal to
\begin{equation}
\label{eq3.1}
\begin{array}{l}
{\ds \sum_{1\le i\le m} \frac{1}{1-q(i)} - \sum_{1\le i<j\le m} \frac{1}{1-q(i,j)}
+ \sum_{0\le i<j<l\le m} \frac{1}{1-q(i,j,l)} + \ldots}
\\
{\ds \ldots + (-1)^{m-g+1}
\sum_{0\le i_1<i_2<\cdots<i_{m-g}\le m}
\frac{1}{1-q(i_1,\ldots ,i_{m-g})}+}
\\
{\ds
+ \sum_{1\le k\le g}(-1)^{m-g+k+1}
{m \choose m-g+k}}
\end{array}
\end{equation}
\end{proposition}

\noindent
\emph{Proof:}\
Let us define, as before, the set of random variables $V_1,\ldots, V_m$, where
$V_i$ denotes the (random) number of groups of coupons that we need to collect in order
to obtain for the first time the $i$-th coupon.
It is immediate to see that the random variables $V_i$'s have now a geometric law with
parameter $1-q(i)$.
Similarly, the random variables $\min(V_i,V_j)$ have a geometric law
with parameter $1-q(i,j)$ and so on up to the random variables
$\min(V_{i_1},\ldots,V_{i_{m-g}})$, which
have a geometric law
with parameter $1-q(i_1,\ldots ,i_{m-g})$.
As in the uniform case,
the minimum of more random variables $V_i$, i.e. $\min(V_{i_1},\ldots,V_{i_{k}})$ for $k>m-g+1$,
will be equal to the constant random variable $1$.
Applying the Maximum-Minimums principle, we obtain that the expected number of
groups of coupons
that we need to complete the collection is equal to
$\EE[\max(V_1, \ldots, V_m)]$, which is equal to (\ref{eq3.1}).
\qed

Let us now assume that a group of coupons could contain
more copies of the same type.
Defining now $\Omega=\{1,\ldots,m\}^g$ and by
$S(i_1,\ldots,i_k)=\{(a_1,\ldots,a_g)\in \Omega:
a_j\not\in \{i_1, \ldots,i_k\} \ \mbox{for} \ j=1, \ldots, g\}$,
we will denote
by $q_\omega, \omega \in \Omega$ the probability
to purchase (at any given time)
the $\omega$ group of coupons.
As before, given $k\in \{1,\ldots,m\}$, we shall denote by
$q(i_1,\ldots,i_k)$ the probability to purchase a group of coupons which
does not contain any of the coupons $i_1,\ldots,i_k$.
As pointed out before, the assignment of the probabilities $q_\omega$ is
in general not simple and most of all not unique.
However, if we assume
to draw without replacement $g$ elements
from a population composed by $m$ different types of individuals
which are present in different proportions,
it is easy to compute the previous probabilities.
To fix the notation, let $N$ be the total number of individuals and
$N_1, \ldots, N_m$ the number of individuals of any given type.
A simple computation gives, for example,
\begin{equation}
\label{ppp1}
q(1) = \prod_{j=0}^{g-1} \frac{N-N_1-j}
{N-j} = \frac{P(N-N_1,g)}{P(N,g)}
\end{equation}
where $P(n,k)$ denotes the number of ordered sequences
of $k$ elements from $n$,
and similarly, fixed $i_1\neq i_2\neq\ldots \neq i_k$,
\begin{equation}
\label{ppp2}
q(i_1,\ldots,i_k) = \prod_{j=0}^{g-1} \frac{N-N_{i_1}-\ldots -N_{i_k}-j}
{N-j}
= \frac{P(N-N_{i_1}-\ldots N_{i_k},g)}{P(N,g)} \ .
\end{equation}

Following the same ideas as before, it
is easy to extend the result of Proposition \ref{3.1}
to the present case:
\begin{proposition}
\label{3.11}
The expected number of groups of coupons that we need to complete
the collection, in the case of unequal probabilities $q_\omega$,
is equal to
\begin{equation}
\label{eq3.2}
\begin{array}{l}
{\ds \sum_{1\le i\le m} \frac{1}{1-q(i)} - \sum_{1\le i<j\le m} \frac{1}{1-q(i,j)}
+ \sum_{0\le i<j<l\le m} \frac{1}{1-q(i,j,l)} + \ldots}
\\
{\ds \ldots + (-1)^{m+1} \frac{1}{1-q(1,\ldots ,m)}}
\end{array}
\end{equation}
\end{proposition}

\begin{corollary}
Let a population of size $N$ be composed of $m$ types of different individuals in
given proportions $N_1/N, \ldots, N_m/N$.
The expected number of independent drawn
without replacement of $g$ individual
that we have to perform in order
to observe at least once any type of individuals
is equal to
\begin{equation}
\label{zzz}
\begin{array}{l}
{\ds \sum_{1\le i\le m} \frac{1}{1- \frac{P(N-N_{i},g)}{P(N,g)}}
- \sum_{1\le i<j\le m} \frac{1}{1-
\frac{P(N-N_{i}-N_j,g)}{P(N,g)}
} + }
\\
{\ds
+ \sum_{0\le i<j<l\le m} \frac{1}{1-
\frac{P(N-N_{i}-N_j-N_k,g)}{P(N,g)}
} + \ldots
 + (-1)^{m+1} \frac{1}{1-
\frac{P(N-N_{1}\ldots - N_m,g)}{P(N,g)}}
}
\end{array}
\end{equation}
\end{corollary}

\

Let us now evaluate (\ref{zzz}) for some specific choices 
of the relative distribution in the population and the
size $g$ of the single sample.
First of all it is important to note that if one 
type of individuals is vary rare with respect to the
others, the value (\ref{zzz}) is very close to the 
expected number of samples that we have to draw in
order to obtain one element of this type.
For example, if we choose
$m=4$ and $N_1=10, N_2=100, N_3= 500, N_4=1000$,
we get that the expected number of independent drawn of two
individuals that we need in order to observe at least once
any of the four types, is approximatively equal to $81.5$,
while the expected number of independent drawn of two
individuals that we need in order to observe one individual 
of the first type is equal to $80.7$.

\begin{figure}[h]
	\includegraphics[width=0.7\textwidth]{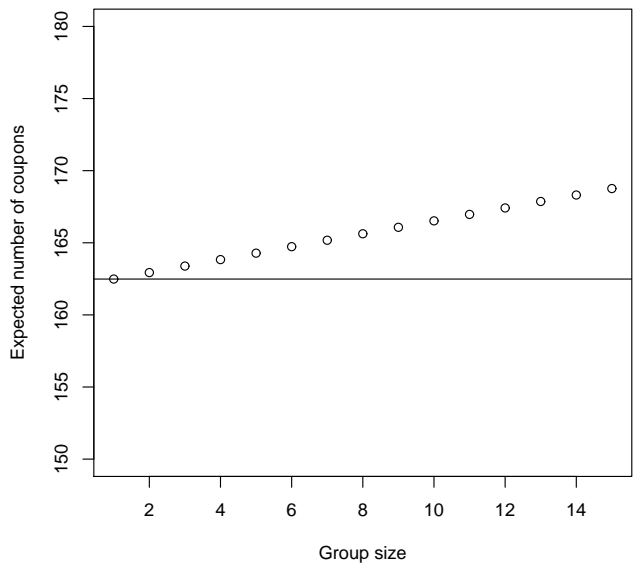}
	\caption{Expected number of individuals to observe at least one individual of any type for different group sizes ($g=1,\ldots,15$), computed using (\ref{zzz}) with $m=4$ and $N_1=10, N_2=100, N_3= 500, N_4=1000$. Comparison with the single arrivals case (solid line)}
	\label{coupon}
\end{figure}

In order to see how the quantity (\ref{zzz}) depends on the size $g$,
we choose again $m=4$ and $N_1=10, N_2=100, N_3= 500, N_4=1000$ and we compute the value of (\ref{zzz}) for $g=1,\ldots,15$. In Figure (\ref{coupon}) we plot the expected numbers of individuals that we have to draw in order to observe at least one individual of any type. As one could expect, this expectation increases with $g$. We also compare these values with the case of single arrivals (solid line).

On the converse, fixed $g=2$, we can consider a population with an increasing
number $m$ of different types. 
%Formula (\ref{zzz}) can be computed just for small value of $m$ (for $m\ge 20$ we are already close to the computability bound); so, for bigger values of $m$, we can approximate the exact value performing a simulation.

\begin{figure}[h]
	\includegraphics[width=0.7\textwidth]{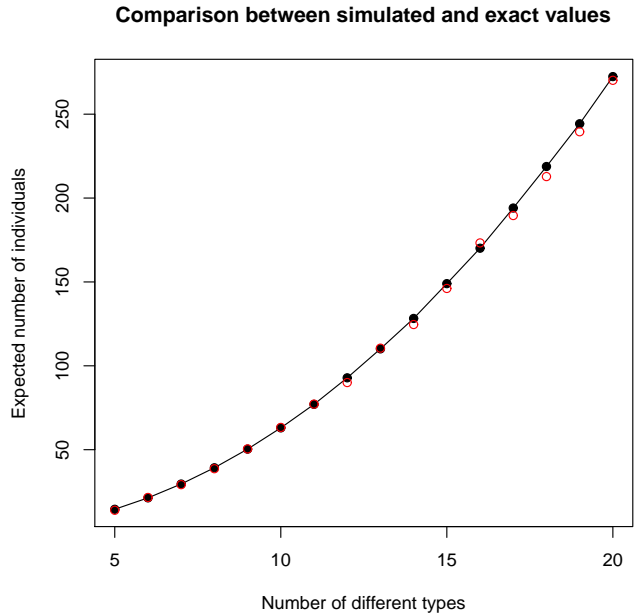}
	\caption{Expected number of independent drawn of two individuals that we need in order to observe at least once
any of the $m$ types, with $m=5,\ldots,20$ (solid black circle). Computation is performed using (\ref{zzz}) with proportion of the types in the population closed to a
Mandelbrot distribution of parameters $c=0.30$ and $\theta=1.75$. Comparison with the simulated values (filled red circle)}
	\label{comparison}
\end{figure}

Taking the proportion of the types in the population closed to a
Mandelbrot\footnote{The Mandelbrot distribution assumes events to be ranked according to their frequency of usage. The $i$-th most probable event has probability $p_i\propto (c+i)^{-\theta}$ for some constant $c\geq 0$ and $\theta$ ranges over $[1,2]$}
distribution of parameters $c=0.30$ and $\theta=1.75$. Figure (\ref{comparison}) shows the exact and the simulated values of (\ref{zzz})
for increasing values of $m$.

\begin{remark}
\label{remark223}
It is important to note that both the expressions
(\ref{eq3.1}) and (\ref{eq3.2}) are computationally hard and
the explicit computation of their values possible just for small values
of $m$.
\end{remark}

%MR1978107 Ross
%MR1082197 Stadje
%MR0061772 Von Schelling
%MR0268929 Rosen
%MR0266273 Rosen asymptotics
%MR959649 Holst

%
\bibliography{coupon} % Here the name of the file .bib
                       % (bibliography database)
                        % that must be used.
\bibliographystyle{plain}

\end{document}